\documentclass[11pt]{article}

\usepackage{shapepar,bm}
\usepackage[dvips,arrow,matrix,ps,color,line]{xy}
\usepackage{theorem,amsfonts,amssymb,amscd}
\usepackage{geometry}
\usepackage{makeidx}

\usepackage{graphicx}
\usepackage{graphics}

\usepackage{epstopdf}

\usepackage{amssymb, graphics,hyperref}

\hypersetup{
    colorlinks = true,
    linkcolor = blue,
    anchorcolor = red,
   citecolor = blue,
    filecolor = red,
    pagecolor = red,
    urlcolor = blue}

\newcommand{\mathsym}[1]{{}}

\usepackage[usenames]{color}
\definecolor{MyLightMagenta}{cmyk}{0.1,0.8,0,0.1}
\definecolor{MyDarkBlue}{rgb}{0.1,0,0.3}

\makeindex

\def\bfh{{\mathbf h}}
\def\bfu{{\mathbf u}}
\def\bfv{{\mathbf v}}

\def\Qcal{\mathcal Q}

\def\Wcal{{\mathcal W}}

\def\ZZ{\mathbb Z}

\def\CC{\mathbb C}

\def\Scal{{\mathcal S}}

\def\QQ{\mathbb Q}
\def\PP{\mathbb P}

\def\cocoa{{\hbox{\rm C\kern-.13em o\kern-.07em C\kern-.13em o\kern-.15em A}}}

\def\bfu{{\bf u}}

\def\blamb{{\bm \lambda}}
\def\bmu{{\bm \mu}}

\def\ep{{\epsilon}}

\def\w2M{\bigwedge^2M}

\def\w{\wedge }

\protect
\protect
\protect
\protect

\protect
\protect

\def\sra{\rightarrow}
\def\lra{\longrightarrow}

\def\proof{\noindent{\bf Proof.}\,\,}
\def\qed{{\hfill\vrule height4pt width4pt depth0pt}\medskip}
\def\be{\begin{equation}}
\def\ee{\end{equation}}
\def\bclm{\begin{claim}}
\def\eclm{\end{claim}}
\def\beqn{\begin{eqnarray}}
\def\eeqn{\end{eqnarray}}
\def\beqn*{\begin{eqnarray*}}
\def\eeqn*{\end{eqnarray*}}

\theoremstyle{change}
\theorembodyfont{\rmfamily}

\newtheorem{claim}{}[section]

\global
\global

\def\no@breaks#1{{\def\\{ \ignorespaces}#1}}    

  \makeatletter
\def\cleardoublepage{\clearpage\if@twoside \ifodd\c@page\else
\hbox{} \thispagestyle{empty}
\newpage
\if@twocolumn\hbox{}\newpage\fi\fi\fi} \makeatother

\usepackage{eso-pic,graphicx}
\makeatletter
\newcommand\BackgroundPicture[2]{%
  \setlength{\unitlength}{1pt}%
  default \put(0,\strip@pt\paperheight){%
  \parbox[t][\paperheight]{\paperwidth}{%
    \vfill
     \centering \includegraphics[angle=#2, width=15cm, height=15cm,  bb=0 0 150 150]{#1}
    \vfill
}}} %
\makeatother

\usepackage{amssymb}
\usepackage{palatino,url}

\title{{\Large "On one Property of one Solution of one Equation"\\
or\\ Linear ODEs, Wronskians and  Schubert Calculus}\thanks{\noindent The work was partially sponsored by PRIN
``Geometria sulle Variet\`a Algebriche" (Coordinatore A.~Verra) and Politecnico di Torino. The second author was
sponsored by an INDAM-GNSAGA grant (2009)  for Visiting Professors at the Politecnico di Torino.}}

\author{ L.~Gatto \and  I.~Scherbak}

\date{}

\begin{document}

\maketitle

\hfill{\it To the Memory of Vladimir Arnold}

\medskip
 \abstract{For a linear ODE with indeterminate coefficients, we exhibit a fundamental system of solutions
 explicitly, in terms of the coefficients. We show that the generalized Wronskians of the
fundamental system are given by an action of the Schur functions on the usual Wronskian, and thence
enjoy Pieri's and Giambelli's formulae. As an outcome, we obtain a natural isomorphism between the free
module generated by the generalized Wronskians and the singular homology module of the Grassmannian.}

\bigskip \hfill "{\it The existence of mysterious relations between all these different domains }

 \hfill {\it is the most striking and delightful feature of mathematics.}"  V.~Arnold,~\cite{ArnFroPer}

\section*{Introduction}

In this paper we deal with the {\it universal linear ODE}
\be u^{(r+1)}-e_1u^{(r)}+\ldots+(-1)^{r+1}e_{r+1}u\,=\,0\,\label{eq:univ} \ee
where $e_1,\ldots , e_{r+1}$ are indeterminates.

Let $E_r=\QQ[e_1,\ldots , e_{r+1}]$ be a  polynomial $\QQ$-algebra in the ODE coefficients, and
$E_r[[t]]$ the $E_r$-algebra of formal power series of $t$, supplied with the standard formal derivation
with respect to $t$. In Sec.~\ref{sec:ude}, we solve equation (\ref{eq:univ}) explicitly in $E_r[[t]]$,
as well as a non-homogenous one.

In  Sec.~\ref{Giambelli}, we use the obtained {\em universal} fundamental system of solutions to express
the {\em generalized Wronskians}, which are numerated by partitions, in terms of the ordinary Wronskian
of (\ref{eq:univ}). It appears, for a given partition the ratio of the corresponding generalized
Wronskian with the usual one is a specialization of the Schur function associated with the same
partition (Theorem~\ref{giambthmwr}). Thus the product in the free module spanned by the generalized
Wronskians obeys the same rule as the product of Schur functions, and hence as the intersection of
Schubert cells. In particular, this relation imply a natural isomorphism between the free
$\ZZ[e_1,\ldots,e_{r{ +1}}]$-module generated by the generalized Wronskians and the
$H^*(G(r,\PP^\infty))$-module of the singular homology of the Grassmannian $G(r,\PP^\infty)$ of
$r$-dimensional linear subvarieties of the infinite-dimensional complex projective space.

\medskip\noindent
The title ``On one Property of one Solution of one Equation'' was mentioned many times by V.~Arnold as
an example of a meaningless one. Therefore one of the authors, the one who was fortunate to be Arnold's
student and to attend Arnold Seminar, dreamed to publish one day a (preferably interesting) paper with
this title, just to surprise Arnold.

In fact, there was one opportunity, when that author, in collaboration with another former student of
V.~Arnold, prepared for publication a paper dedicated to Arnold's 65-th birthday. But the collaborator
was strict and rejected this idea immediately.

This time the dream nearly comes true, but not really, as no chance to surprise Arnold anymore...

\claim{\bf  Universal linear ODE.}\ \ In Sec.~\ref{sec:ude}, we revise the theory of linear ODE's in an
algebraic context.  For the classical theory see e.g.~\cite{arnold}.

\medskip\noindent
Denote by $D:E_r[[t]]\sra E_r[[t]]$ the standard formal derivative with respect to $t$. Throughout the paper we
use the following notation: \be f=f(t)=\sum_{n\geq 0}{\displaystyle a_n\over n!}t^n:=\sum_{n\geq 0}a_n{t^n\over
n!}\,, \ \ D^if=f^{(i)}(t):=\sum_{n\geq 0}a_{n+i}{t^n\over n!}\,,\ \ f^{(i)}(0):=a_i\,{ ,} \label{eq:fmlpwser}
\ee for $f\in E_r[[t]]$, $a_n\in E_r$.

\medskip\noindent
Our {\tt one Equation} is (\ref{eq:univ}).

\medskip\noindent
Define $h_j\in E_r$ \ ($j\in\ZZ$)\  by means of the following generating function:
\be
{1\over 1-e_1t+e_2t^2-\ldots +(-1)^{r+1}e_{r+1}t^{r+1}}\,=\,\sum_{j\in\ZZ}\, h_jt^j\,. \label{eq:e-h} \ee In
particular, $h_j=0$ for $j<0$, $h_0=1$, $h_1=e_1$, $h_2=e_1^2-e_2$, etc.

\medskip\noindent
The origin of the present work is the following observation:

\begin{center}{\it  The formal power series\ $\displaystyle{u_0:=\sum_{n\geq 0} h_n{t^n\over n!}}$\
is a solution of {\rm(\ref{eq:univ})}.}
\end{center}
This is our {\tt  one Solution}.

\medskip\noindent
For each  $1\leq j\leq r$, define the formal power series $u_j=u_j(t)$  as the unique element in the ideal
$(t^j)$ of $E_r[[t]]$ such that $D^ju_j=u_0$. The explicit formulae are
\be
 u_j\,=\,\sum_{n\geq j}\,h_{n-j}\,{\displaystyle{t^n\over n!}}\,,\qquad 0\leq j\leq r \,. \label{eq:fundsys}
\ee

\medskip
\noindent {\bf Theorem~A.}\ \ {\em The power series {\rm(\ref{eq:fundsys})} form a fundamental system of
solutions of  {\rm(\ref{eq:univ})}.}

\medskip\noindent
 This is our {\tt one Property}.

\medskip\noindent
 We call  (\ref{eq:fundsys})
 the {\it universal} fundamental system. The reason is the following.
 Denote by $E_r[T]$ the $E_r$-polynomial ring in the indeterminate $T$. Equation (\ref{eq:univ}) may be
written then as $U_{r+1}(D)u=0$, where
$$U_{r+1}(T)\,=\,T^{r+1}\,-\,e_1T^r\,+\,\ldots\, +\, (-1)^{r+1}e_{r+1}$$
is the {\it universal} polynomial of the degree $r+1$.
 For any $\QQ$-algebra $A$ and for any monic polynomial $P(T)\in A[T]$, we can solve the ODE
$P(D)u=0$  in  formal power series $A[[t]]$.

\medskip
\noindent {\bf Theorem~B.}\ \ {\em If $P(T)\in A[T]$ is a monic polynomial of degree $r+1$, then there is a
unique $\QQ$-algebra homomorphism $\psi:E_r\sra A$  mapping  the universal fundamental
system~{\rm(\ref{eq:fundsys})} to a fundamental system of the linear ODE $P(D)v=0$.}

\medskip\noindent
Our version of the theory of linear homogeneous and non-homogeneous ODE's is described in Sec.~\ref{sec:ude}.

\medskip
\claim{\bf Derivatives of the Wronskian and Generalized Wronskians.}\ \ For a non-decreasing set of $r+1$
integers, $\lambda_0\geq \lambda_1\geq \ldots\geq \lambda_r\geq 0\,, $  denote $\blamb=(\lambda_0,\lambda_1,
\ldots,\lambda_r)$ the partition of $|{\blamb}|:=\lambda_0+\lambda_1+\ldots + \lambda_r$.

\medskip\noindent Given an $(r+1)$-tuple ${\mathbf f}$  of formal power series,
$${\mathbf f}:=\pmatrix{f_0\cr f_1\cr \vdots\cr f_r}, $$
and a partition $\blamb$, the {\it generalized Wronskian} $W_\blamb({\mathbf f})$ is defined as the determinant
of the matrix whose $i$-th row entries are the derivatives of order $j+\lambda_{r-j}$ of $f_i$, for \ $0\leq
i,j\leq r$. The usual Wronskian $W({\mathbf f})=W_0({\mathbf f})$ corresponds to the trivial partition
$\lambda_j=0\,, 0\leq j\leq r$,
 \be
 W_\blamb({\mathbf f})\,:=\,\det
(f_i^{(j+\lambda_{r-j})})_{0\leq i,j\leq r}\,,\ \ W_0({\mathbf f})\,:=\,\det (f_i^{(j)})_{0\leq i,j\leq r}\,.
\label{eq:ww} \ee

\medskip\noindent
It is not surprising that the derivatives of all orders of the Wronskian are positive $\ZZ$-linear combinations
of generalized Wronskians.  However the intriguing fact was  that {\it the coefficients of the linear
combinations have a precise interpretation in terms of enumerative geometry of subspaces of the complex
projective space}.

\medskip\noindent
For example, the dimension of  the Grassmann variety $G(r, \PP^d)$ of  $r$-dimensional subspaces in  the complex
projective $d$-dimensional space is $(r+1)(d-r)$. It was observed in~\cite{Gat2}, that the coefficient of
$W_{(d-r, d-r,\ldots,d-r)}({\mathbf f})$ in the $\ZZ$-linear combination of  the $(r+1)(d-r)$-th derivative of
$W({\mathbf f})$ is the Pl\"ucker degree of $G(r, \PP^d)$ (cf.~\cite[Example 14-7-11(iii)]{Fu1}).

\medskip\noindent
Here we put such a remark in a proper framework. Recall that partitions are described by means of Young--Ferrers
diagrams, and a standard Young tableau is a filling in the Young--Ferrers diagram of $\blamb$ with numbers
$1,\ldots, |\blamb|$ in such a way that the numbers in any column and in any row increase, \cite{Fu2}.

\medskip\noindent{\bf Proposition.} \ \ {\em  We have
$$ D^{k+1}W(\mathbf f)\,=\, \sum_{|\blamb|=k} c_\blamb W_\blamb(\mathbf f)\,,
$$
where $c_\blamb$ is the number of the standard Young tableaux of the Young--Ferrers diagram $\blamb$. \qed}

\medskip\noindent Numbers $c_\blamb$'s and their interpretation in terms of Schubert calculus (\ref{clambda})
are very well known, see e.g., \cite{Fu2}.  In particular, these numbers can be calculated by the {\it hook
formula},
$$
c_\blamb\,=\,{|\lambda|! \over k_1\cdot\,\ldots\, \cdot k_{|\lambda|}}\,,
$$
where $k_j$'s are the {\it hook lengths} of the boxes of $\blamb$.

\claim {\bf Wronski--Schubert Calculus.}\ \ In Sec.~\ref{Giambelli}, the $(r+1)$-tuple ${\mathbf f}$  is a
fundamental system  of solutions of the ODE (\ref{eq:univ}), i.e., a basis of $\ker U_{r+1}(D)$.  One can
re-write the equation (\ref{eq:univ}) in the form
\[
\left|\matrix{u & u'&\ldots & u^{(r+1)}\cr f_0 & f_0' & \ldots & f_0^{(r+1)}\cr f_1 & f_1' & \ldots &
f_1^{(r+1)}\cr \vdots & \vdots & \ldots & \vdots \cr  f_r & f_r' & \ldots & f_r^{(r+1)}}\right|=\,0\,.
\]
In notation (\ref{eq:ww}), we get
\[
W_0(\mathbf f)u^{(r+1)}-W_{(1)}(\mathbf f)u^{(r)}+W_{(1^2)}(\mathbf f)u^{(r-1)}+\, \ldots\, +
(-1)^{r+1}W_{(1^{r+1})}(\mathbf f)u\,=\,0\,,
\]
where  $(1^k):=(\underbrace{1,\ldots, 1}_{k\,\,\mathrm{times}},0\ldots,0)$ is the {\it primitive} partition of
$k$.

\noindent Thus we have $W_{(1^k)}(\mathbf f)=e_kW(\mathbf f)$,\ $1\leq k\leq r+1$. In particular,
$W_{(1)}(\mathbf f)$ is the derivative of $W_0(\mathbf f)$, i.e., the Wronskian solves the first order universal
equation $u'-e_1u=0$. This is the well-known Liouville theorem, cf. ~\cite[Ch.~3, \S 27.6]{arnold}.

\medskip\noindent
 Relation (\ref{eq:e-h}) is satisfied, in particular, by
 the {\it complete symmetric functions} $h_j$'s and the {\em elementary symmetric functions}
 $e_j$'s  of $r+1$ variables (see e.g., \cite[Ch.~I, \S 2]{MacDonald}). As the Referee pointed out,
 this interpretation of $h_j$'s and $e_j$'s fits the approach of A.~Lascoux,
\cite{lascoux}, within symmetric functions are treated as operators in polynomial spaces.

Being motivated by Schubert Calculus problems, we have chosen to interpret  (\ref{eq:e-h}) in terms of
characteristic classes, as the relation between Chern classes (our $e_j$'s) and Segre classes (our
$h_j$'s) of the tautological bundle over the Grassmannian $G(r,\PP^\infty)$. We refer to  \cite{Fu1} for
the intersection theory of complex Grassmannians.

\medskip\noindent
It is well-known that $G=G(r,\PP^\infty)$ possesses a cellular decomposition into  Schubert cells
parameterized by the partitions $\blamb=(\lambda_0,\lambda_1, \ldots, \lambda_r)$. The homology classes
of their closures, $\{\Omega_\blamb\}$, form a $\ZZ$-basis of $H_*(G, \ZZ)$.
 The $\cap$ product map turns  the homology $H_*(G,\ZZ)$ into a { free module of rank $1$} over the cohomology,
 generated by the fundamental class $[G]$. The Schubert classes  $\sigma_\blamb\in H^*(G,\ZZ)$ are then defined
 via the equality $\Omega_\blamb=\sigma_\blamb\cap [G]$, for every partition $\blamb$.  In $H^*(G,\ZZ)$ the
multiplication (or $\cup$ product, or intersection) of cohomology classes is defined, and the classical
{\em Pieri's formula} gives the product of a {\it special} Schubert class $\sigma_{k}$, i.e., that
corresponding to the {\it special} partition $(k,0,\ldots,0)$\,  of $k$,  with an arbitrary one,
\[
\sigma_{k}\sigma_{\blamb}=\sum_\bmu \sigma_\bmu\,,\ \ 1\leq k\leq
r+1,
\]
where the sum is taken over all the partitions $\bmu=(\mu_0,\mu_1, \ldots, \mu_r)$ such that
\be
|\bmu|=k+|\blamb|\,,\ \ \mu_0\geq\lambda_0\geq\mu_1\geq\lambda_1\geq\ldots\geq \mu_r\geq \lambda_r\,.
\label{eq:pieri}
\ee
Notice that  $c_{\blamb}$'s from the Proposition above are coefficients in the decomposition of $\sigma_1^k$
into the sum of Schubert classes:
\be
\sigma_1^k\,=\,\sum_{|\blamb|=k} c_\blamb \sigma_\blamb\,.\label{clambda}
\ee

\medskip\noindent
The special Schubert classes ${\bm\sigma}: =(\sigma_0, \sigma_1,\ldots, \sigma_k, \ldots)$ are multiplicative
generators of the cohomology ring  $H^*(G,\ZZ)$. Recall that the Schur polynomial  corresponding to the
partition $\blamb$ is:
\be
 \Delta_\blamb({\bm x }):= \det(x_{\lambda_{r-j}+j-i})_{0\leq i,j\leq r}\,, \label{eq:schur} \ee
where ${\bm x}=(x_k)_{k\in\ZZ}$ are variables, \cite{MacDonald}.

\medskip\noindent
The classical  {\em Giambielli's formula} gives any Schubert class $\sigma_\blamb$ as the specialization of the
corresponding Schur polynomial at the special Schubert classes:
 \be
\sigma_\blamb=\Delta_\blamb({\bm\sigma})\,, \label{eq:delta}
\ee
where  one sets $\sigma_j:=0$ for $j<0$.

\medskip\noindent
In Sec.~\ref{Giambelli} we prove the  Giambielli's  and the Pieri's formulae for generalized Wronskians
(Theorem~\ref{giambthmwr}, Corollary~\ref{piericor}).  We write ${\bf h}:=(h_j)_{j\in\ZZ}$ for the sequence
$h_j$'s defined in~(\ref{eq:e-h}).

\medskip
\noindent {\bf Theorem~C.}\ \  {\em We have
$$W_\blamb(\mathbf f)=\Delta_\blamb(\bfh)W_0(\mathbf f)\,,\ \
h_kW_\blamb(\mathbf f)=\sum_\bmu W_\bmu(\mathbf f)\,, \ k\geq 1\,,
$$
where the sum is over all the partitions $\bmu$ satisfying {\rm (\ref{eq:pieri})}. }

\medskip\noindent
Denote by  $\Wcal(\mathbf f)$ the free $\ZZ$-module generated by the generalized Wronskians $W_\blamb(\mathbf
f)$'s.

\medskip
\noindent {\bf Corollary.}  \ \ {\em The correspondence $\Omega_\blamb\mapsto W_\blamb(\mathbf f)$ defines an
isomorphism between $\Wcal(\mathbf f)$  as a $H^*(G,\ZZ)\cong \ZZ[e_1,\ldots, e_r]$-module and the
$H^*(G,\ZZ)$-module $H_*(G,\ZZ)$, the singular homology of the infinite Grassmannian}.

 \claim{\bf Comments.}\ \  Wronskians are ubiquitous in mathematics. We especialy are interested in their
role in algebraic geometry. For example, any linear system on the projective line defines an
$(r+1)$-dimensional subspace of  the vector space of complex polynomials in one indeterminate, i.e. an
element of the corresponding Grassmannian. The {\em Wronski map} sending the polynomial subspace to the
class modulo $\CC^*$ of the Wronski determinant of its basis appears. In particular, the critical points
of a rational function in one variable are roots of the Wronskian of the denominator
 and the numerator, the ramification points of linear systems
on curves occur as the zero locus of a certain Wronskian section of suitable line bundles, etc. In this context,
the Wronski map has been studied by many authors, see \cite{EH, EreGab1, sottile, Sch}.

The moduli points corresponding to curves of fixed genus possessing special Weierstrass points provide
another example; they occur as the zero locus of a suitable derivative of a Wronskian associated to the
sheaf of relative differentials, \cite{GatScand, GatPon}. Further relations between Wronskians of linear
systems and Schubert calculus on a Grassmann bundle are investigated in~\cite{CuGaNi, EreGab2,
GatSalehyan}.

Discovering of Wronskians in the Bethe Ansatz of the $sl_n$ Gaudin model was crucial both for study of
Bethe vectors and for pure algebro-geometric problems, like transversality of the intersection of
Schubert varieties, \cite{MTV, MV, Scherb1, Scherb3}. In~\cite{Sch}, relation to the $sl_2$
representation theory was used to calculate the intersection number of some Schubert varieties in the
Grassmannian  of projective lines; later on  a purely algebro-geometric proof, based on the Wronskian
inspired methods of ~\cite{GatSant}, appeared in ~\cite{NBFSC}.

Generalized Wronskians   have also appeared within different contexts,  e.g. in the theory of
Weierstrass points on curves (\cite{Schmidt, towse}), and in connection with number theory
(\cite{MMO,lacunary}).

Recently, D.~Laksov and A.~Thorup showed that the $n$-th exterior power of a polynomial ring in $n$
indeterminates is a free module over its subring of symmetric functions, and that the module structure is
equivalent to the Schubert Calculus for infinite Grassmannians, \cite{LakTh, LakTh1}. Our construction may be
considered as a concrete realization of this equivalency.

\section{Solution to Linear ODEs in { Formal} Power Series.}\label{sec:ude}

First of all we solve the universal equation (\ref{eq:univ}) with coefficients in $E_r$, i.e., find $\ker
U_{r+1}(D)\subset E_r[[t]]$. Substitution of (\ref{eq:fmlpwser}) into (\ref{eq:univ}) gives the following
condition on the coefficients of $t^n$.

\bclm{\bf Proposition.} \label{condsolpr} {\em The series
$f=\sum_{n\geq 0}\, x_n\displaystyle{t^n\over n! }\,\in\ker U_{r+1}(D)\,$
 if and only if}
\be x_{n+1}-e_1x_{n}+\ldots+(-1)^{r+1}e_{r+1}x_{n-r}=0\,, \ \
n\geq r\,.\label{eq:condsol}  \ee \qed \eclm
This recurrence relation on the coefficients $x_n$'s imposes no
restrictions on $x_0,x_1,\ldots, x_r$; they are the {\em initial conditions} of the solution.

\bclm{\bf Theorem.}\label{base} {\em The $E_r$-module $\ker U_{r+1}(D)\subset E_r[[t]]$  is free of rank
$r+1$ and generated by $u_0,\ldots ,u_r$ given by {\rm(\ref{eq:e-h})}, {\rm(\ref{eq:fundsys})}. }
\eclm

\proof First we will prove that every $u_j$ is a solution of (\ref{eq:univ}).

\medskip\noindent
 Let us begin with  $u_0$. The coefficients of $t^n$ on  both sides of (\ref{eq:e-h}) are
the same, that is
\be
h_0=1\,,\ \ h_{n+1}-e_1h_{n}+\ldots+(-1)^{r+1}e_{r+1}h_{n-r}=0\,,\ n\geq 0\,,\label{eq:recurs}
\ee
and Proposition~\ref{condsolpr} immediately implies $u_0\in\ker U_{r+1}(D)$.

Next, let us show  $u_r\in \ker U_{r+1}(D)$. Indeed, $u_{r-j}=D^ju_r$ \ $(0\leq j\leq r)$, and we get
\[
D^{r+1}u_r-e_1\cdot D^ru_r+\ldots+(-1)^{r+1}e_{r+1}u_r=
\]
\[
=u_0'-e_1u_0+e_2u_{1}+\ldots+(-1)^{r+1}e_{r+1}u_r=
\]
\[
=\sum_{n\geq 0} (h_{n+1}-e_1h_{n}+e_2h_{n-1}+\ldots+(-1)^{r+1}e_{r+1}h_{n-r}){t^n\over n!}.
\]
By~(\ref{eq:recurs}), all the coefficients of the expansion  vanish.

Finally for each $0\leq i\leq r$,
\[
U_{r+1}(D)u_i=U_{r+1}(D)D^{r-i}u_r=D^{r-i}U_{r+1}(D)u_0=0,
\]
proving that  $u_i$ is a solution as well.

\medskip\noindent
Now we will show that $u_0,\ldots ,u_r$ generate $\ker U_{r+1}(D)$ as an $E_r$-module.
Assume that  $f=\sum_{n\geq
0}\,x_n\displaystyle{t^n\over n!}$ is a solution of $U_{r+1}(D)u=0$.

Define $\Lambda_0,\ldots,\Lambda_r\in E_r$ as the unique solution of the linear system \be
\pmatrix{1&0&\ldots&0\cr h_1&1&\ldots&0\cr \vdots&\vdots&\ddots&\vdots\cr h_r&h_{r-1}&\ldots&1}\,
\pmatrix{\Lambda_0\cr\Lambda_1\cr\vdots\cr \Lambda_r}=\pmatrix{x_0\cr x_1\cr\vdots\cr x_r}\,.
\label{eq:lambda-x} \ee We contend that \be
f\,=\,\Lambda_0u_0+\Lambda_1u_{1}+\ldots+\Lambda_ru_r\,,\label{eq:flincomb}
\ee or, in terms of coefficients, \be
x_{n}=\Lambda_0h_n+\Lambda_1h_{n-1}+\ldots+\Lambda_rh_{n-r}\,,\ \ n\geq r+1\,. \label{eq:dncl} \ee The proof is
by induction on $n\geq r+1$.  For $n=r+1$ we have, according to (\ref{eq:condsol}),
\[
x_{r+1}\,=\,e_1x_r-\ldots+(-1)^re_{r+1}x_0\,.
\]
First we express $x_{r+1}$ via $\Lambda_j$'s from (\ref{eq:lambda-x}), and then apply (\ref{eq:recurs}):
\begin{eqnarray*}
 x_{r+1}&=&e_1(h_r\Lambda_0+\ldots +h_1\Lambda_{r-1}+\Lambda_r)-\ldots+(-1)^re_{r+1}\Lambda_0=\\
&=&\Lambda_0(e_1h_r+e_2h_{r-1}-\ldots +(-1)^re_{r+1})+\ldots+\Lambda_re_1=\\
&=&\Lambda_0h_{r+1}+\ldots+\Lambda_rh_1\,.
\end{eqnarray*}
Suppose now (\ref{eq:recurs}) true for all $r\leq m\leq n$. Then
\begin{eqnarray*}
x_{n+1}&=&e_1x_n-e_2x_{n-1}+\ldots-(-1)^{r+1}e_rx_{n-r}=\\
&=&e_1(\sum_{i=0}^r\Lambda_ih_{n-i})-e_2(\sum_{i=0}^r\Lambda_ih_{n-1-i})+\ldots -(-1)^{r+1}e_{r+1}
(\sum_{i=0}^r\Lambda_ih_{n-r-i})=\\
&=&\sum_{i=0}^r\Lambda_i(e_1h_{n-i}-e_2h_{n-1-i}+\ldots -(-1)^{r+1}e_{r+1}h_{n-r-i})=
\\
&=&\sum_{i=0}^r\Lambda_ih_{n+1-i}=\Lambda_0h_{n+1}+\Lambda_1h_{n}+\ldots+ \Lambda_rh_{n-r+1},
\end{eqnarray*}
as desired.

Finally, we conclude that $u_0,\ldots, u_r$ are linearly independent. Indeed, for the trivial solution $x_j=0$
for all $j\geq 0$; hence  $\Lambda_0=\ldots =\Lambda_r=0$ is the unique solution of the homogeneous linear
system (\ref{eq:lambda-x}). \qed

The inverse matrix of the matrix in (\ref{eq:lambda-x}) is well-known (see e.g.~\cite[Ch.~I, \S 2]{MacDonald}),
so we can explicitly find $\Lambda_j$'s in terms of $x_j$'s:
\[
\pmatrix{\Lambda_0\cr\Lambda_1\cr\Lambda_2\cr\vdots\cr \Lambda_r}=\pmatrix{1&0&0&\ldots&0\cr
-e_1&1&0&\ldots&0\cr e_2&-e_1&1&\ldots&0\cr\vdots&\vdots&\vdots&\ddots\cr
(-1)^re_{r}&(-1)^{r-1}e_{r-1}&(-1)^{r-2}e_{r-2}&\ldots&1}\pmatrix{x_0\cr x_1\cr x_2\cr\vdots\cr x_r}\,.
\]
Substitution into (\ref{eq:flincomb}) gives

\bclm{\bf Corollary (Universal Solution of the Cauchy Problem).} \label{cauchythm} \ \ {\em
Let $y_0,y_1,\ldots, y_r$ be indeterminates over  $E_r$. The unique  solution
of~(\ref{eq:univ}) over $E_r[y_0,y_1,\ldots,y_r]$ satisfying  $D^if(0)=y_i$, $0\leq i\leq r$, is as follows:}
\be
\matrix{g(t)=\Lambda_0(y)u_0(t)+\Lambda_1(y)u_1(t)+\ldots +\Lambda_r(y)u_r(t)\,, \cr\cr
\Lambda_j(y)=y_j-e_1y_{j-1}+\ldots +(-1)^je_jy_0\,, \ \ \ 0\leq j\leq r\,.}\label{eq:univcauchy}
\ee\qed
\eclm

\claim{\bf Universality.}\label{psi}
 For any $\QQ$-algebra $A$,  denote by $A[T]$ the $A$-algebra of polynomials of $T$ and by $A[[t]]$  the
$A$-algebra of formal power series of $t$. For $a_1,\ldots , a_{r+1}\in A$,  a $\QQ$-algebra homomorphism $\psi:
E_r\sra A$ mapping $e_j\mapsto a_j$ $(1\leq j\leq r+1)$ naturally induces homomorphisms of $\QQ$-algebras $
\widehat{\psi}: E_r{ [}T{ ]}\sra A{ [}T{ ]}$ and $\widetilde{\psi}: E_r[[t]]\sra A[[t]]$. Now we are in position
to solve the ODE $P(D)v=0$, for any monic polynomial
\[
P(T)=T^{r+1}-a_1T^r+\ldots+(-1)^{r+1}a_{r+1}\in A[T]\,.
\]
Indeed, let $\psi$ be  the unique $\QQ$-algebra homomorphism $E_r\sra A$ sending $e_i\mapsto a_i$, \ $1\leq
i\leq r+1$. Then $P(T)$ is the image of $U_{r+1}(T)$ under $ \widehat{\psi}$. Clearly, $\widetilde{\psi}$ maps
any solution of the universal ODE $U_{r+1}(D)u=0$ to a solution of $P(D)v=0$. Moreover, as the matrix  of
(\ref{eq:lambda-x}) is unimodular, its image is unimodular as well, and so the proof of Theorem~\ref{base} may
be repeated verbatim for the images of $u_j(t)$'s. We get

\bclm{\bf Theorem.}\label{univsol} {\em  (1) The series $f(t)=\sum_{n\geq 0}\,x_n\displaystyle{t^n\over n!}\,
\in A[[t]]$
is a solution of $P(D)v=0$ if and only if
\[
x_{n+1}-a_1x_{n}+\ldots+(-1)^{r+1}a_{r+1}x_{n-r}=0\,,\ \  n\geq r+1\,.
\]

\noindent(2)\ $\ker P(D)$ is a free $A$-module of rank $r+1$
generated by  $v_0,\ldots ,v_r$, the image of
{\rm(\ref{eq:fundsys})}, via the unique homomorphism of
$\QQ$-algebras $E_r\sra A$ mapping $e_i\mapsto a_i$\,, $1\leq
i\leq r+1$. In particular, the fundamental system $v_j$'s
satisfies $D^jv_r=v_{r-j}$, $0\leq j\leq r$.

\noindent(3)  The unique solution of $P(D)v=0$  with given initial
conditions $x_0,\ldots, x_r\in A$ is \be
f(t)=V_0(x)v_0(t)+V_1(x)v_1(t)+\ldots
+V_r(x)v_r(t)\,,\label{eq:specialization} \ee where  $V_j(x)$'s
are the images of $\Lambda_j(y)$'s under the unique $\QQ$-algebra
homomorphism
$$
E_r[y_0,y_1,\ldots, y_r]\sra A
$$
mapping $e_i\mapsto a_i$ and $y_j\mapsto x_j$. \qed }
\end{claim}

\bigskip Consider as an example the equation $v''-3v'+2v=0$\,
(clearly the real and the complex numbers are $\QQ$-algebras). The
corresponding (characteristic) polynomial is
$T^2-3T+2=(T-2)(T-1)$, and the standard prescription,
\cite{arnold}, gives the fundamental system $ f_1=e^t\,,
f_2=e^{2t}$. According to Theorem~\ref{univsol}, $f_1=v_0-2v_1$
and $f_2=v_0-v_1$.

Of course, one can easily obtain this relation directly. Indeed,
the substitution $e_1=3\,, e_2=2$ into (\ref{eq:e-h}) gives the
image of the universal fundamental system,
$$v_0\,=\,1+3t+{7t^2 \over 2}+{15t^3\over 3!}+\dots\,,\ \
v_1\,=\,t+{3t^2\over 2}+{7t^3 \over 3!}+{15t^4\over 4!}+\dots\,,
$$
i.e., $v_0=2f_2-f_1\,, v_1=f_2-f_1$ and  $v_1'(t)=v_0(t)$.

\claim{\bf Universal Euler formula.}\label{eul}\ \  For $x\in A$, define $\exp(xt):=\sum_{n\geq 0}\,
x^n{t^n\over n!}$. If $x$ is a root of $P(T)$, it is not  surprising  that $\exp(xt)\in \ker P(D)$. Indeed, the
condition of  Theorem~\ref{univsol}~(1) holds:
\[
x^{n+1}-a_1x^{n}+\ldots+(-1)^{r+1}a_{r+1}x^{n-r}\,=\,x^{n-r}P(x)\,=\,0\,, \ \  n\geq r+1\,.
\]
In particular, if $\epsilon\in\CC$ is a primitive root of $(-1)$ of order $r+1$, then  $\exp(\epsilon t)$ is the
solution of the ODE $u^{(r+1)}+u=0$ over $\CC$ with initial conditions $1,\epsilon,\ldots,\epsilon^r$. According
to  Theorem~\ref{univsol}~(3),
\be
\matrix{\exp(\epsilon\,t)=v_0(t)+\epsilon v_{1}(t)+\ldots +\epsilon^rv_r(t)\,, \cr\cr
 v_j(t)=\sum_{k=0}^\infty\, (-1)^k{\displaystyle {t^{j+k(r+1)} \over (j+kr+k)!}}\,,\ \  \  0\leq j\leq r\,.} \label{eq:eulerouniv}
\ee
For $r=1$, the image of the universal fundamental system $(u_0,u_1)$ is $(cos\,t, sin\,t)$, $\epsilon=i$, and
the classical  Euler formula $\exp(it)=\cos\,t + i\sin\,t$ follows.

More generally, let $\alpha:=T+(U_{r+1}(T))$ be the {\em universal
root} of $U_{r+1}(T)$, cf.~\cite{Lak}. Define
$E_r[\alpha]:=E_r[T]/(U_{r+1}(T))$. The polynomial $U_{r+1}(T)$ is
defined over $E_r[\alpha]$ as well, and $U_{r+1}(\alpha)=0$. In
fact $E_r[\alpha]$ is the {\em universal splitting algebra} of
$U_{r+1}(T)$ as a product of two monic polynomials, one of degree
$1$. It is universal in the sense that for any $E_r$-algebra $B$
where $U_{r+1}(T)$ splits as the product $(T-\beta )q(T)$, with
$\beta\in B$ and $q(T)\in B[T]$ monic of degree $r$, there is a
unique homomorphism $\phi: E_r[\alpha]\sra B$ mapping
$\alpha\mapsto \beta$, \cite{Lak}. Thus $exp(\alpha t)$ is a
solution of $U_{r+1}(D)u=0$ defined over $E_r[\alpha]$, and we
have \be \matrix{\exp(\alpha
t)=u_0(t)+\Lambda_1(\alpha)u_1(t)+\ldots+\Lambda_r(\alpha)u_r(t)\,,\cr\cr
\Lambda_j(\alpha)=\alpha^j-e_1\alpha^{j-1}+\ldots +(-1)^je_j\,,\ \
0\leq j\leq r\,.} \label{eq:univeuler} \ee Again, for $r=1$, one
has $\exp(\alpha t)=u_0(t)+(\alpha-e_1)u_1(t)\,$. Under the unique
$\QQ$-algebra homomorphism $E_1[\alpha]\sra \CC$ defined by
$e_1\mapsto 0$, $e_2\mapsto 1$ and $\alpha\mapsto i$, this formula
becomes the classical Euler formula $\exp(it)=\cos\,t + i\sin\,t$.

\claim {\bf Remark on fundamental systems.}\label{rem:fs}\ \ In the algebraic context under consideration, there
is a difference between a fundamental system, which is a basis of $\ker P(D)$, and a set $f_0,\ldots , f_r$ of
linearly independent solutions of $P(D)u=0$.

\medskip
For example, consider $A=\QQ[a_1,a_2]$ and the equation $u''-(a_1+a_2)u'+a_1a_2u=0$. The characteristic
polynomial is $P(T)=(T-a_1)(T-a_2)$, and hence $\exp{a_1t}$, $\exp{a_2t}$ are two linearly independent solutions
of the equation. However these solutions do not form a basis of  $\ker P(D)$ over $A$. Indeed, consider $v_0\,,
v_1$, the image of the universal fundamental system { $(u_0,u_1)$ of $U_1(D)u=0$} via  the homomorphism $E_2\sra
A$, sending $e_1\mapsto a_1+a_2$ and $e_2\mapsto a_1a_2$. By  Theorem~\ref{univsol}, $(v_0\,, v_1)$ is a basis
of  $\ker P(D)$ over $A$. The relation to  $\exp{a_1t}$, $\exp{a_2t}$ is given by formula~(\ref{eq:univeuler}):
\[
\exp{a_1t}=v_0-a_2v_1\,,\ \ \exp{a_2t}=v_0-a_1v_1.
\]
 Thus $\exp{a_1t}-\exp{a_2t}=(a_1-a_2)v_1$, but $(a_1-a_2)$ is not invertible in $A$! We conclude that there is
no way to  express $v_1$  as an $A$-linear combination of $\exp{a_1t}$ { and} $\exp{a_2t}$.

However, for the same equation {\it considered over the localization} $B:=A\left[1/(a_1-a_2)\right]\,$  both
$(v_0, v_1)$  and $(\exp{a_1t}$, $\exp{a_2t})$  are  fundamental systems of $U_2(D)y=0$ over $B$, i.e., bases of
$\ker P(D)$ over $B$. Since $(v_0,v_1)$ is also a fundamental system  over $A$, it may be considered as "the
most economical" one, see~\cite{GatScherb} for further discussion.

\medskip
Recall that for any  algebra $A$ we have (\cite{cartan}):

\medskip
\centerline{\em $g=g(t)\in A[[t]]$ is invertible in $A[[t]]$ if and only if $g(0)$ is invertible in $A$.}

\medskip\noindent Let $P(T)\in A{ [}T{ ]}$ be a monic polynomial of degree $r+1$ and $f_0,\ldots, f_r\in \ker P(D)$.
Denote by $C:=(f^{(j)}_i(0))_{0\leq i<j\leq r}$  {\em the initial condition matrix} of $f_i$'s. Then $\det C\in
A$ is the constant term of the determinant of the transition matrix  from  the image of the universal
fundamental system, $v_i$'s, to $f_j$'s.

\bclm{\bf Proposition.}\label{basis}\ \  {\em (1) \ The solutions $f_0,\ldots, f_r$ of  $P(D)y=0$
form a basis of $\ker P(D)$ over $A$ if and only if $\det C$ is invertible in $A$.

\noindent (2) If $\det C\neq 0$ and non-invertible in $A$, then $f_0,\ldots, f_r$ are $A$-linearly independent,
but do not form a fundamental system. They do form a fundamental system over  any $\QQ$-algebra extension of
$B:=A[1/\det C]$.} \qed
\eclm
In particular, if $P(T)$ has $(r+1)$ distinct roots $a_0,\ldots, a_r$,  then $\exp(a_0t),\ldots, \exp(a_rt)$
form a fundamental system only if the Vandermonde determinant
$$\prod_{0\leq i<j\leq r}(a_i-a_j)$$
is invertible in $A$. On the other hand, the image of the universal fundamental system is always a fundamental
system.

\claim{\bf Non-homogeneous case.}\ \ Here we discuss the case of the non-homogeneous ODE. Let
\be
f=\sum_{n\geq 0}a_n{t^n\over n!}\in E_{r+1}[[t]]\,.
\ee
Consider the Cauchy problem:
\be
\matrix{U_{r+1}(D)y=f\,,\cr\cr D^{k}y(0)=b_k{\rm \ \ for}\,\, k=0,1,\ldots, r.} \label{eq:nonheq}
\ee

\bclm{\bf Proposition.}\label{nonh} {\em\  Let $\sum_{n\geq 0}p_{n+r+1}t^n$ be the formal power series defined
by the following generation function:
\[
{\sum_{n\geq 0}a_nt^n\over (1-e_1t+\ldots+(-1)^{r+1}e_{r+1}t^{r+1}}\,=\,\sum_{n\geq 0}p_{n+r+1}t^n\,.
\]
Then
\[
u_p=\Lambda_0(b)u_0(t)+\Lambda_1(b)u_1(t)+\ldots +\Lambda_r(b)u_r(t)\,+\,\sum_{n\geq r+1}p_n{t^n\over n!}\,,
\]
where $\Lambda_j(b)=b_j-e_1b_{j-1}+\ldots+(-1)^je_jb_0$ for $0\leq j\leq r$, is the unique solution of the
Cauchy problem~{\rm(\ref{eq:nonheq})}. }
\eclm

 \proof\ \
Let us collect coefficients of $t^n$ on  both sides of the given equality,
\[(p_{r+1}+p_{r+2}t+p_{r+3}t^2+\ldots)(1-e_1t+\ldots+(-1)^{r+1}e_{r+1}t^{r+1})=a_0+a_1t+a_2t^2+\ldots\,.
\]
Similarly to the homogeneous case (see Theorem~\ref{base}), we get
\be
p_{r+1+n}-p_{r+n}e_1+\ldots+(-1)^{r+1}p_ne_{r+1}=a_n\,, \ \ n\geq 0,\label{eq:secondenh}
\ee
and this exactly means that $y_0(t)=\sum_{n\geq r+1}p_n{t^n\over n!}$ is a solution of $U_{r+1}(D)y=f$ with
initial conditions $y_0^{(k)}(0)=0$ for $k=0,1,\ldots,r$.

Therefore the solution of the Cauchy problem will be the sum of $y_0$ and the solution $y_b$ of
 $U_{r+1}(D)y=0$ that satisfies the initial conditions $y_b^{(k)}(0)=b_k$ for $k=0,1,\ldots, r$.
 By Corollary~\ref{cauchythm}, one can take  $y_b(t)=\Lambda_0(b)u_0(t)+\Lambda_1(b)u_1(t)+\ldots +\Lambda_r(b)u_r(t)$.

The solution $u_p=y_0+y_b$ is clearly unique, as the difference of two solutions of the same Cauchy problem is
the (unique) solution of the homogeneous equation that satisfies the zero initial conditions. \qed

\section{\bf Generalized Wronskians and  Schubert Calculus }\label{Giambelli}

Results of this section hold for any fundamental system of solutions of the universal ODE $U_{r+1}(D)u=0$. They
are formulated in terms of generalized Wronskians, and  therefore it is enough to prove them for one fundamental
system. To ease computations, the choice of the  universal fundamental system is the more convenient one.

 From now on $\bfu=\bfu(t)$ will denote the ``column vector'' of the universal solutions (\ref{eq:fundsys})
of $U_{r+1}(D)u=0$,
\begin{center}
$ \bfu(t)=\pmatrix{u_0\cr u_1\cr\vdots\cr u_r}. $
\end{center}

For any partition $\blamb=(\lambda_0,\lambda_1, \ldots, \lambda_r)$, the generalized  Wronskian
$W_{\blamb}(\bfu)$ defined by (\ref{eq:ww}), has the form
\[
W_{\blamb}(\bfu)=D^{\lambda_r}\bfu\w D^{1+\lambda_{r-1}}\bfu\w\ldots\w D^{r+\lambda_{0}}\bfu\,.
\]
We use notation (\ref{eq:schur}) and write ${\bfh}=(h_k)_{k\in\ZZ}$ for  $h_k$'s  given by (\ref{eq:e-h}).

\bclm{\bf Proposition.}\label{lambwron} { \em We have
\be
\matrix{W_\blamb(\bfu)=\displaystyle{\sum_{n\geq 0}\sum_{|{\bmu}|=n}{n\choose
\bmu}\Delta_{\blamb+\bmu}(\bfh){t^n\over n!}}\,,\cr\cr W_\blamb(\bfu)(0)=\Delta_\blamb(\bfh)\,,}
\label{eq:genwrpart}
\ee
where ${\bmu}=(\mu_0, \mu_1,\ldots, \mu_r)$ is a partition of $|{\bmu}|=\mu_0+ \mu_1+\ldots +\mu_r$ and
\[
{n\choose \bmu}={n\choose \mu_0,\mu_1,\ldots,\mu_r}={n!\over \mu_0!\mu_1!\cdot \ldots\cdot \mu_r!}\,.
\]
}
\eclm

\proof   By definition of ${\bfu}$, see (\ref{eq:fundsys}),
 \[
 \bfu=\sum_{n\geq 0}\pmatrix{h_n\cr h_{n-1}\cr\vdots\cr h_{n-r}}\cdot{t^n\over  n!}\,,
 \]
{ and then}  the $(i+\lambda_{r-i})$-th derivative is
\[
 \displaystyle{D^{i+\lambda_{r-i}}\bfu=\sum_{\mu_{r-i}\geq 0}\pmatrix{h_{\lambda_{r-i}+\mu_{r-i}+i}\cr
h_{\lambda_{r-i}+\mu_{r-i}+i-1}\cr\vdots\cr h_{\lambda_{r-i}+\mu_{r-i}+i-r}}{t^{\mu_{r-i}}\over \mu_{r-i}!}}\,,
\]
where we put $\mu_{r-i}$ instead of $n$. Thus
\[
W_\blamb(\bfu)=D^{\lambda_{r}}\bfu\w D^{1+\lambda_{r-1}}\bfu\w\ldots\w D^{r+\lambda_0}\bfu=
\]
\[
=\sum_{\bmu}\pmatrix{h_{\lambda_{r}+\mu_{r}}\cr h_{\lambda_{r}+\mu_{r}-1}\cr\vdots\cr
h_{\lambda_{r}+\mu_{r}-r}}{t^{\mu_r}\over \mu_r!}\w\pmatrix{h_{\lambda_{r-1}+\mu_{r-1}+1}\cr
h_{\lambda_{r-1}+\mu_{r-1}}\cr\vdots\cr h_{\lambda_{r-1}+\mu_{r-1}-(r-1)}}{t^{\mu_{r-1}}\over
\mu_{r-1}!}\w\ldots\w \pmatrix{h_{\lambda_{0}+\mu_{0}+r}\cr h_{\lambda_{0}+\mu_{0}-1}\cr\vdots\cr
h_{\lambda_{0}+\mu_{0}}}{t^{\mu_{0}}\over \mu_{0}!}=
\]
\[
=\sum_{\bmu}\pmatrix{h_{\lambda_{r}+\mu_{r}}\cr h_{\lambda_{r}+\mu_{r}-1}\cr\vdots\cr
h_{\lambda_{r}+\mu_{r}-r}}\w \pmatrix{h_{\lambda_{r-1}+\mu_{r-1}+1}\cr h_{\lambda_{r-1}+\mu_{r-1}}\cr\vdots\cr
h_{\lambda_{r-1}+\mu_{r-1}-(r-1)}}\w\ldots\w \pmatrix{h_{\lambda_{0}+\mu_{0}+r}\cr
h_{\lambda_{0}+\mu_{0}+r-1}\cr\vdots\cr h_{\lambda_{0}+\mu_{0}}}\cdot{t^{\mu_r}\over \mu_r!}{t^{\mu_{r-1}}\over
\mu_{r-1}!}\cdot\ldots\cdot {t^{\mu_{0}}\over \mu_{0}!}=
\]
\[
=\sum_{n\geq 0}\,\,\,\sum_{\{\bmu: |\bmu|=n\}}{n!\over \mu_0!\mu_1!\ldots\mu_r!}\left|\matrix{h_{\lambda_{r}+\mu_{r}}&
h_{\lambda_{r-1}+\mu_{r-1}+1}&\ldots&h_{\lambda_{0}+\mu_{0}+r} \cr
h_{\lambda_{r}+\mu_{r}-1}&h_{\lambda_{r-1}+\mu_{r-1}}&\ldots&h_{\lambda_{0}+\mu_{0}+r-1}\cr\vdots&\vdots&\ddots&\vdots\cr
h_{\lambda_{r}+\mu_{r}-r}&h_{\lambda_{r-1}+\mu_{r-1}-(r-1)}&\ldots& h_{\lambda_{0}+\mu_{0}}}\right| {t^n\over
n!}=
\]
\[
=\sum_{n\geq 0}\,\,\,\sum_{\{|\bmu|=n\}}{n\choose \bmu}\Delta_{\blamb+\bmu}(\bfh){t^n\over n!}
\]
as desired.\qed

\bclm{\bf Theorem.}\label{giambthmwr} {\em Giambelli's formula for generalized Wronskians holds, i.e., for every
partition $\blamb$ we have
\begin{center}
$W_\blamb(\bfu)=\Delta_\blamb({\bfh})W_0(\bfu)\,.$
\end{center}
}
\eclm

\proof Clearly  $W_\blamb(\bfh)$ is an $E_r$-multiple of $W_0(\bfh)$. Indeed,  each column of the form $D^k\bfu$
with $k\geq r+1$  occurring  in the expression for $W_\blamb(\bfu)$   can be replaced with a linear combination
of lower derivatives of $\bfu$,  by~(\ref{eq:univ}) and its consequences. Thus one obtains the product of an
element of $E_r$ with the determinant  involving derivatives of $\bfu$ of orders at most $r$ only. This
determinant  is $W_0({ \bfu})$, up to a permutation of the columns. Hence $W_\blamb(\bfu)=\gamma_\blamb
W_0(\bfu)$ for some $\gamma_\blamb\in E_r$. But the ratio of the constant terms of  $W_\blamb({ \bfu})$ and
$W_0({ \bfu})$, according to   Proposition~\ref{lambwron}, is
\[
\gamma_\blamb={ {W_{\blamb}(\bfu)(0)\over W_{0}(\bfu)(0)}}=\Delta_\blamb({\bfh})\,,
\]
 and the statement follows.
 \qed

\bclm{\bf Corollary.}\label{piericor} {\em Pieri's formula for generalized Wronskians holds, i.e., for every
partition $\blamb$ and $k\geq 0$ we have
\begin{center}
$ h_kW_\blamb(\bfu)=\sum_\bmu W_\bmu(\bfu)\,, $
\end{center}
where the sum is taken over all the partitions $\bmu$ satisfying~{\rm(\ref{eq:pieri})}.  }
\eclm

\proof By ~\cite[Lemma A.9.4]{Fu1}, we have $h_k\Delta_\blamb(\bfh)=\sum_\bmu\Delta_\bmu(\bfh)$, where the sum
is over all the partitions $\bmu$  satisfying (\ref{eq:pieri}). Applying Theorem~\ref{giambthmwr}, we get
\begin{center}
$ h_kW_\blamb(\bfu)=(h_k\Delta_\blamb(\bfh))W_0(\bfu)=\sum_\bmu\Delta_\bmu(\bfh)W_0(\bfu)=\sum_\bmu
W_\bmu(\bfu)\,, $
\end{center}
as desired.\qed

\claim{\bf Schubert Calculus} (\cite{Fu1}, \cite{GH}). \ \

Let   $G:=G(r, \PP^d)$ be the complex Grassmann variety parameterizing $r$-dimensional subspaces of
$\PP^d$ and let
\begin{center}
$
0\lra \Scal_r\lra G\times \CC^{d+1}\lra \Qcal_r\lra 0
$
\end{center}
be the universal exact sequence over $G$, where $\Scal_r$ is the rank $r+1$ {\em universal subbundle} of
$G\times \CC^{d+1}$ and $\Qcal_r$ the {\em universal quotient bundle} of rank $d-r$.

As it is well known, the ring of symmetric functions in $r+1$ indeterminates surjects onto the
cohomology ring of $G$ by mapping the $i$-th complete polynomial to the $i$-th Chern class of the
universal quotient bundle, and the kernel of the map gives relation between generators, see e.g.
\cite{Fu1}.

In terms of Schubert Calculus on Grassmannians, the integral cohomology of $G$ is a $\ZZ$-algebra
finitely generated by the Chern classes $(-1)^i\ep_i:=c_i(\Scal_r)\in H^*(G,\ZZ)$ of the tautological
bundle or, due to relation $c(\Scal_r)c(\Qcal_r)=1$, by the Chern classes $\sigma_i:=c_i(\Qcal_r)\in
H^*(G,\ZZ)$ of $\Qcal_r$. Clearly $\ep_i=0$ if $i>r+1$ and $\sigma_i=0$ if $i>d-r$. Denote $
{\bm\sigma}_r:=(1,\sigma_1,\sigma_2,\ldots)$, and let $\sigma_\blamb:=\Delta_\blamb({\bm\sigma}_r)\in
H^{2|\blamb|}(G,\ZZ)$.

The Grassmann variety  $G$ possesses a cellular decomposition $\{B_\blamb\}$, with respect to some complete
flag of linear subvarieties of $\PP^d$,  and each affine cell $B_\blamb$ has real codimension $2|\blamb|$. The
singular homology classes $\Omega_\blamb\in H_{2(r+1)(d-r)-2|\blamb|}(G,\ZZ)$ of the closures of $B_\blamb$ in
$G$ form a $\ZZ$-basis of $H_*(G,\ZZ):=\oplus H_{2(r+1)(d-r)-2|\blamb|}(G,\ZZ)$.  The key fact of Schubert
calculus on Grassmannian is that $\sigma_\blamb$ is the Poincar\'e dual of $\Omega_\blamb$:
\[
\sigma_\blamb\cap[G]=\Omega_\blamb\,,
\]
where $[G]$ is the fundamental class of the Grassmann variety, and  the $\cap$ product turns the homology into a
{ free module of rank $1$}  over the cohomology, { generated by the fundamental class}.

\claim{\bf Wronski Calculus.} Consider the differential equation
\be
D^{r+1}v-\ep_1D^rv+\ldots+(-1)^{r+1}\ep_{r+1}v=0\label{eq:UDESCH}
\ee
 over  the $\QQ$-algebra $H^*(G, \QQ)=H^*(G,\ZZ)\otimes_\ZZ\QQ$.
Let $\psi: E_r\sra H^*(G,\QQ)$ denote  the  unique $\QQ$-module
 homomorphism defined by $e_i\mapsto\ep_i$.  By Theorem~\ref{univsol}, the image under the induced homomorphism
 $\widetilde{\psi}: E_r[[t]]\sra H^*(G,\QQ)[[t]]$ of the universal fundamental system
 $\bfu$ is a fundamental system  $\bfv$ of solutions  of~(\ref{eq:UDESCH}):
 $$
 \bfv=\pmatrix{v_0\cr v_1\cr\vdots\cr v_r}\,,  \ \ v_k=\widetilde{\psi}(u_k)\,, \ \ 0\leq k\leq r\,.
 $$
In particular,  $\psi (h_k)= {\sigma}_k$, and  Theorem~\ref{giambthmwr} asserts
\[
W_\blamb(\bfv)=\sigma_\blamb\cdot W_0(\bfv)
\]
or, more suggestively,
\[
\sigma_\blamb={W_\blamb(\bfv)\over W_0(\bfv)}.
\]

\medskip
Define $\Wcal(\bfv):=\bigoplus_{\blamb}W_\blamb(\bfv)$. It  is  a free  $H^*(G,\ZZ)$-module of rank $1$
generated by $W_0(\bfv)$.

\bclm {\bf Proposition.} \ \  {\em The $\ZZ$-module isomorphism
\[
wr:H_*(G,\ZZ)\sra \Wcal_\blamb(\bfv),
\]
 mapping $\Omega_\blamb\mapsto W_\blamb(\bfu)$,  is an isomorphism of $H^*(G,\ZZ)$-modules.}
\eclm

\proof \  It suffices  to show that
 $wr(\sigma_k\cap\Omega_\blamb)=\sigma_kW_\blamb(\bfv)$, for each $k\geq 0$.
According to Corollary~(\ref{piericor}), we have $\sigma_k\sigma_\blamb=\sum_\bmu\sigma_\bmu$,
 where the sum is over all partitions $\bmu$  satisfying~(\ref{eq:pieri}). Therefore
\begin{eqnarray*}
wr(\sigma_k\cap \Omega_\blamb)&=&wr(\sigma_k\cap (\sigma_\blamb\cap [G])=wr ((\sigma_k\cup
\sigma_\blamb)\cap[G])=\\ \hskip90pt  &=&wr(\sum_\bmu\sigma_\bmu\cap[G])=\sum_\bmu W_\bmu(\bfv)=
\sigma_kW_\blamb(\bfv).\hskip90pt \qed
\end{eqnarray*}
Thus the $\cap$ product  can be interpreted as the product of $h_k$ with a generalized Wronskian, and
the $\cup$-product in cohomology as the product in the ring $\ZZ[\ep_1,\ldots, \ep_{r+1}]$. Notice that
that $\ep_i$'s are not necessarily algebraically independent. In fact $\ZZ[\ep_1,\ldots, \ep_{r+1}]$ is
the quotient of $E_r$ through the unique ring epimorphism defined by $e_i\mapsto \ep_i$.

\medskip \noindent In fact, the polynomial $\QQ$-algebra $E_r$ can be interpreted as the cohomology of the
infinite Grassmannian $G(r, \PP^\infty)$, through the telescoping construction, similarly
to~\cite[p.~302]{BottTu}. More precisely, for each $d\geq r$, there is a natural inclusion
 $G(r,\PP^d)\hookrightarrow G(r, \PP^\infty)$. The corresponding arrow reversed $\QQ$-algebra map
$H^*(G(r,\PP^\infty))\sra H^*(G(r, \PP^d)$ is  nothing else but our unique map $\psi$ sending
$e_i\mapsto \ep_i$.  Hence the induced map $\widetilde{\psi}$ maps the universal generalized Wronskians
$W_\blamb(\bfu)$'s to the generalized Wronskians $W_\blamb(\bfv)$'s. Since $G(r, \PP^\infty)$ is an
infinite $CW$-complex, and since the homology $H_*(G(r,\PP^\infty), \ZZ)$ is generated by the classes of
the closures of the Schubert cells, it is thence clear that $\Wcal(\bfu)$ is a module over $H^*(G(r,
\PP^\infty))=\ZZ[e_1,e_2,\ldots, e_{r+1}]$ and is $\ZZ$-isomorphic to $H_*(G(r,\PP^\infty),\ZZ)$. Again,
this extends to an isomorphism of $H^*(G(r,\PP^\infty),\ZZ)$-modules.

\medskip
\noindent {\bf Acknowledgments.} The first author thanks M.~Codegone and C.~Cumino,  the second one thanks
M.~Kazarian, for precious help and discussions. For the second author, it is a pleasure to thank the {
Dipartimento di Matematica del  Politecnico di Torino and INDAM-GNSAGA} for hospitality and nice working
conditions. The authors also thank each other for support and cooperation.

\parbox[t]{3in}{{\rm Letterio~Gatto}\\
{\tt letterio.gatto@polito.it}\\
{\it Dipartimento~di~Matematica}\\
{\it Politecnico di Torino}\\
{\it ITALY}} \hspace{1.5cm}
\parbox[t]{2.5in}{{\rm Inna~Scherbak}\\
{\tt scherbak@post.tau.ac.il}\\
{\it School of Mathematical Sciences}\\
{\it Tel Aviv University}\\
{\it ISRAEL}}

\end{document}